\newfont{\bcb}{msbm10}
\newfont{\matb}{cmbx10}
\newfont{\got}{eufm10}

\documentclass[12pt]{article}
\usepackage[cp1250]{inputenc}

\usepackage{amsmath,amsthm}
\usepackage{amssymb,latexsym}
\usepackage{enumerate}

\newcommand{\matP}{\mathbb{P}}
\newcommand{\matR}{\mathbb{R}}

\newcommand{\matN}{\mathbb{N}}

\begin{document}

\title{On hereditarily rational functions}

\author{by Krzysztof Jan Nowak}




\date{}

\maketitle


\begin{abstract}
In this paper, we give a short proof of a theorem by Koll\'{a}r on
hereditarily rational functions. This is an answer to his appeal
to find an elementary proof which does not rely so much on
resolution of singularities. Our approach does not make use of
desingularization techniques. Instead, we apply a stronger version
of the \L{}ojasiewicz inequality. Moreover, this allows us to
sharpen Koll\'{a}r's theorem.
\end{abstract}


In his recent paper~\cite{Kol}, Koll\'{a}r introduced a class of
continuous rational functions on an algebraic, possibly singular,
variety $X$. Those functions, called hereditarily rational, are
defined by the condition that their restrictions to each algebraic
subvariety $Y$ of $X$ remain rational. He proved that every
continuous rational function on a smooth algebraic variety is
hereditarily rational (Proposition~8). Continuous rational
functions on smooth algebraic varieties were investigated by
Kucharz~\cite{Ku}. Also,
Fichou--Huisman--Mangolte--Monnier~\cite{F-H-M-M} examined
regulous functions on singular algebraic varieties, i.e.\ those
functions which extend to continuous rational functions on an
ambient, smooth algebraic variety. For the rudiments of real
algebraic geometry, we refer the reader to~\cite{BCR}.

The significance of the class of hereditarily rational functions
is visible especially in view of a theorem by
Koll\'{a}r~\cite{Kol} which indicates that hereditarily rational
functions enjoy the good properties of continuous rational
functions on smooth algebraic varieties. Below, we state and prove
a sharpening of Koll\'{a}r's theorem to the effect that one can
find an extension $F$ of a given hereditarily rational function
$f$ with the same indeterminacy locus.

\vspace{2ex}

{\bf Theorem on hereditarily rational functions.}
\begin{em}
Let $X$ be a subvariety of a smooth algebraic variety $M$ and $f:
X \longrightarrow \matR$ a continuous rational function. Then the
two conditions are equivalent:

i) $f$ is hereditarily rational;

ii) $f$ extends to a continuous rational function $F: M
\longrightarrow \matR$.

\noindent Moreover, if $f$ is hereditarily rational and regular
off an algebraic subvariety $Z \subset X$, then we can find its
extension $F$ which is regular on $M \setminus Z$.
\end{em}

\vspace{2ex}

{\bf Remarks.} 1) The additional conclusion about the
indeterminacy locus $Z$ is an essential sharpening of Koll\'{a}r's
theorem, because he constructs an appropriate extension $F$ from a
regular function which descends to $F$ through a finite sequence
of blowings-up biregular over $M \setminus X$. Thus one can only
deduce that $F$ is regular on $M \setminus X$. It seems that the
use of desingularization or transformation to a normal crossing by
blowing up along smooth centers leads to non-trivial modifications
both of the indeterminacy locus of a given rational function $f$
and the singular locus of the algebraic variety $X$.

2) The implication $ii) \Rightarrow i)$ follows immediately from
the above-mentioned proposition. The converse implication remains
valid when $M$ is an arbitrary, possibly singular, algebraic
variety, because it can be embedded into a smooth variety as a
closed subvariety. In particular, it holds if $X$ is an algebraic
subvariety of a Zariski locally closed subset of $\matR^{n}$.

3) Koll\'{a}r's proof of the implication $i) \Rightarrow ii)$
depends heavily on desingularization. He writes that it would be
desirable to find an elementary proof, one that does not rely so
much on resolution of singularities. This article provides a proof
which satisfies, we think, the above demands. We do not make use
of resolution of singularities at all. Instead, we apply a version
of the \L{}ojasiewicz inequality from~\cite{BCR}, Theorem~2.6.6,
recalled below. It also holds for continuous definable functions
in an arbitrary polynomially bounded, o-minimal structure.  Such a
version was formulated and applied in our paper~\cite{Now}, which
is devoted to carrying over the composite function theorem to the
quasianalytic settings.

\vspace{2ex}

{\bf \L{}ojasiewicz Inequality.}
\begin{em}
Let $f,g: A \longrightarrow \matR$ be two continuous
semi-algebraic functions on a locally closed, semi-algebraic subset $A$ of
$\matR^{n}$ such that
$$ \{ x \in A: f(x)=0 \} \subset \{ x \in A: g(x)=0 \}. $$
Then there exist a positive integer $k$ and a continuous semi-algebraic function
$h: A \longrightarrow \matR$ such that $g^{k} = fh$.
\end{em}

\vspace{2ex}

Now, we can readily prove the theorem. Clearly, $f$ is
hereditarily rational iff there exists a filtration
$$ \emptyset = X_{0} \subset X_{1} \ldots X_{m-1} \subset X_{m} =
   X, $$
where $X_i{}$ are algebraic subvarieties of $X$ such that $X_{i}$
is nowhere dense in $X_{i+1}$ and the restriction of $f$ to
$X_{i+1} \setminus X_{i}$ is a regular function for each
$i=0,1,\ldots,m-1$. The proof is by induction with respect to the
dimension of the variety $X$. It is evident if $X$ is of dimension
$0$. The induction step comes down to the following

\vspace{2ex}

{\bf Lemma.}
\begin{em}
Consider two algebraic subvarieties $A \subset X$ of $M$ and a
continuous rational function $f: X \longrightarrow \matR$ regular
on $X \setminus A$ and vanishing on $A$. Then $f$ extends to a
continuous rational function $F: M \longrightarrow \matR$ regular
on $M \setminus A$.
\end{em}

\vspace{2ex}

Obviously, $f$ can be presented as a fraction $g/h$, where $g,h$
are regular functions on $X$, $h \geq 0$ and
$$ \{ x \in X: h(x)=0 \} \subset A. $$
One can find, of course, their regular extensions
$G,H: M \longrightarrow \matR$ such that $H \geq 0$ and
$$ \{ x \in M: H(x) = 0 \} = \{ x \in A: h(x)=0 \} \subset
   A. $$
The rational function $G/H$ is no longer continuous in general.

Consider the blowing-up $\sigma: \widetilde{M} \longrightarrow M$
with respect to the ideal $(G,H)$. Then the pull-back
$$ F_{1} := \frac{G^{\sigma}}{H^{\sigma}} : \widetilde{M} \longrightarrow \matP_{1}
   \ \ \mbox{ with } \
   G^{\sigma} := G \circ \sigma, \ H^{\sigma} := H \circ \sigma, $$
is a regular mapping into the projective line. Let $Y$ be the
birational transform of $X$, $B:= \sigma^{-1}(A)$ and $C$ be the
Euclidean closure of $Y \setminus B$; $C$ is, of course, a closed
semialgebraic subset of $\widetilde{M}$. Clearly, $\sigma$ is a
biregular mapping of $\widetilde{M} \setminus B$ onto $M \setminus
A$.

Observe further that, in the vicinity of $C$, $F_{1}$ is a regular
function leading into $\matR$. Indeed, for a point $ b \in C
\setminus B$, the denominator $H^{\sigma}(b) \neq 0$ whence
$F_{1}(b) \neq \infty$. On the other hand, if $b \in C \cap B$, a
sequence of points $b_{\nu} \in Y \setminus B$, $\nu \in \matN$,
tends to $b$. But then the sequence $a_{\nu} := \sigma(b_{\nu})
\in X \setminus A$, $\nu \in \matN$, tends to $a := \sigma(b) \in
A$. Consequently,
$$ F_{1}(b) = \lim_{\nu \rightarrow \infty} \, F_{1}(b_{\nu}) =
   \lim_{\nu \rightarrow \infty} \, \frac{G^{\sigma}(b_{\nu})}{H^{\sigma}(b_{\nu})}
   = \lim_{\nu \rightarrow \infty} \, \frac{G(a_{\nu})}{H(a_{\nu})} = $$
$$ = \lim_{\nu \rightarrow \infty} \, \frac{g(a_{\nu})}{h(a_{\nu})}
   = \lim_{\nu \rightarrow \infty} \, f(a_{\nu})  = f(a) = 0 \neq \infty. $$
Hence the pole set $F_{1}^{-1}(\infty)$ of the mapping $F_{1}$ is
disjoint with $C$, as asserted.

Now, take regular functions $P,Q: \widetilde{M} \longrightarrow
\matR$ such that $P,Q \geq 0$ and
$$ B = \{ x \in \widetilde{M}: P(x)=0 \} \ \ \mbox{ and } \ \
   Y = \{ x \in \widetilde{M}: Q(x)=0 \}. $$
Since
$$ \{ x \in \widetilde{M}: H^{\sigma}(x) \} \subset B  \ \ \mbox{ and } \ \
   Y \cap (\widetilde{M} \setminus C) \subset B \cap (\widetilde{M} \setminus C), $$
it follows from the above version of the \L{}ojasiewicz inequality
that there exists a positive integer $k \in \matN$ such that
$$ P^{k} \leq H^{\sigma} \cdot \mbox{\rm constant} \ \ \ \mbox{and}
   \ \ \ P^{k} \leq Q \cdot \mbox{\rm constant} $$
locally on $\widetilde{M} \setminus C$; it means that the
constants in the above inequalities depend on a neighbourhood of a
given point from $\widetilde{M} \setminus C$. Put
$$ F_{2} := \frac{P^{2k}}{P^{2k} + Q} \cdot F_{1}. $$
Observe that the first factor takes on value $1$ on $Y \setminus
B$, whence
$$ F_{2}(y) = (f \circ \sigma)(y) \ \ \mbox{ for all } \  y \in Y
   \setminus B. $$

On $\widetilde{M} \setminus C$, we get
$$ F_{2} = \frac{P^{k}}{P^{2k} + Q} \cdot \frac{P^{k}}{H^{\sigma}}
   \cdot G^{\sigma}, $$
with the first two factors locally bounded on $\widetilde{M}
\setminus C$ and regular off $B$. Therefore, the restriction of
$F_{2}$ to $\widetilde{M} \setminus C$ extends continuously
through $B$ by taking on zero value, because the third factor
$G^{\sigma}$ vanishes on $B$. On the other hand, in the vicinity
of $C$, $F_{2}$ is the product of the regular function $F_{1}$
vanishing on $B$ and a bounded (by $1$) rational function which is
regular off $B$. Consequently, in the vicinity of $C$, $F_{2}$
extends continuously through $B$ by taking on zero value.

Summing up, we see that $F_{2}: \widetilde{M} \longrightarrow
\matR$ is a continuous rational function regular off $B$ and
vanishing on $B$. Consequently, $F_{2}$ is constant on the fibres
of the proper mapping $\sigma$. Hence it descends to a continuous
function $F:\matR^{n} \longrightarrow \matR$ which vanishes on $A$
and is regular off $A$, $F_{2} = F^{\sigma}$. The last assertion
holds since $\sigma$ is biregular over $\matR^{n} \setminus A$. We
have already seen that the functions $F_{2}$ and $f \circ \sigma$
coincide on $Y \setminus B$. Therefore the functions $F$ and $f$
coincide on $X \setminus A$, and thus $F$ is an extension of $f$
we are looking for. Hence the lemma follows and the proof is
complete.

\vspace{2ex}

\begin{small}
\begin{sc}
Institute of Mathematics

Faculty of Mathematics and Computer Science

Jagiellonian University

ul.~Profesora \L{}ojasiewicza 6

30-348 Krak\'{o}w, Poland

{\em e-mail address: nowak@im.uj.edu.pl}
\end{sc}
\end{small}


\vspace{2ex}


\begin{thebibliography}{9}

\bibitem{BCR}
J.~Bochnak, M.~Coste, M.-F.~Roy, {\em Real Algebraic Geometry\/},
Ergebnisse der Mathematic und ihrer Grenzgebiete, Vol.~36,
Springer Verlag, Berlin, 1998.

\bibitem{F-H-M-M}
G.~Fichou, J.~Huisman, F.~Mangolte, J.-P.~Monnier, {\em Fonctions
r\'{e}gulues\/}, arXiv 1112.3800 [math.AG].

\bibitem{Kol}
J.~Koll\'{a}r, {\em Continuous rational functions on real and
p-adic varieties\/}, arXive 1101.3737 [math.AG].

\bibitem{Ku}
W.~Kucharz, {\em Rational maps in real algebraic geometry\/},
Adv.\ Geom.\ {\bf 9} (4) (2009), 517--539.

\bibitem{Now}
K.J.~Nowak, {\em A note on Bierstone--Milman-Paw\l{}ucki's paper
''Composite differentiable functions''\/}, Ann.\ Polon.\ Math.\
{\bf 102} (3) (2011), 293--299.

\end{thebibliography}
\end{document}